\documentclass[12pt]{spieman}  % 12pt font required by SPIE;
\usepackage{amsmath,amsfonts,amssymb}
\usepackage{graphicx}
\usepackage{setspace}
\usepackage{tocloft}
\usepackage{changepage}
\usepackage{booktabs}
\usepackage{soul, color, xcolor}

\usepackage{lineno}
\usepackage[colorlinks=true]{hyperref}
\hypersetup{urlcolor=blue, citecolor=blue}
\usepackage{hyperref}
\allowdisplaybreaks
\usepackage{algorithm}
\usepackage{algorithmic}
\title{Weighted structure tensor total variation for image denoising}

\author[a]{Xiuhan Sheng}
\author[a]{Lijuan Yang}
\author[a,*]{and Jingya Chang}
\affil[a]{Guangdong University of Technology, School of mathematics and statistics, 161 Yinglong Road, Tianhe District, Guangzhou, China, 510520}

\cftpagenumbersoff{figure}
\cftpagenumbersoff{table}

\begin{document}
%\linenumbers
\maketitle

\begin{abstract}
For image denoising problems, the structure tensor total variation (STV)-based models  show good performances when compared with other competing regularization approaches. However, the STV regularizer  does not  couple the local information of the image and may not maintain the  image details.  Therefore, we employ the anisotropic weighted matrix  introduced in the anisotropic total variation (ATV) model to improve the STV model. By applying the weighted matrix to the discrete gradient of the patch-based Jacobian operator in STV, our proposed weighted STV (WSTV) model can effectively capture local information from images and maintain their details during the denoising process. The optimization problem in the model is solved by a fast first-order gradient projection algorithm with a complexity result of $O(1 / i^2)$. For images with different Gaussian noise levels, the experimental results demonstrate that the WSTV model can effectively improve the quality of restored images compared to other TV and STV-based models.
\end{abstract}

% Include a list of up to six keywords after the abstract

% Include email contact information for corresponding author
{\noindent \footnotesize\textbf{*}Address all correspondence to Jingya Chang,  \linkable{jychang@gdut.edu.cn} }

\begin{spacing}{2}   % use double spacing for rest of manuscript

\section{Introduction}
\label{sect:intro}  % \label{} allows reference to this section

As an important area in digital image processing\cite{1}, image restoration  covers topics of image denoising\cite{2,3,4}, deblurring\cite{5,6,7}, medical imaging\cite{8,9,10}, super-resolution reconstruction\cite{11,12,13}, and so on. During the processes of acquisition and transmission, noise signals inevitably contaminate digital images. The purpose of image restoration is to restore clear images from degraded images and simultaneously  preserve  image details, such as textures and edges. In this paper, we mainly focus on image denoising.

The variational method is one of the most widely used strategies for  image denoising problems. For models in the variational method, a key point is the regularization term. Among the existing regularization terms, total variation (TV)\cite{14} is the most common one. TV can maintain sharp edges during noise removal, which overcomes the disadvantage of excessively smooth edges caused by Tikhonov regularization\cite{15}. The central disadvantage of  TV is that it may bring artifact boundaries in smooth regions\cite{16}. Numerous models are proposed to deal with such kinds of staircase effects; for instance,
the high-order total variation (HOTV)\cite{17,18}, the total generalized variation  (TGV)\cite{19,20}, and the nolocal total variation (NLTV)\cite{21}. Compared to TV, these modified models exhibit better restoration results and reduce staircase effects  in some regions of the image.

In addition, Leflammatis et al. proposed  the structure tensor total variation (STV) regularization family\cite{22}, which is a modified version of  TV  as well. The difference between the STV model and the abovementioned TV-based models is that STV proposes to penalize the $\ell_{p}$ norms of the square root of the eigenvalues of the structure tensor  rather than the gradient magnitude.  The STV model provides a more robust and richer measure of image variation than the TV-based methods do by exploiting additional information from the neighborhood of every point. However, the isotropic STV model treats the  discrete gradients in the horizontal and vertical directions equally, which ignores the variation in local features of the image. Ref.~\citenum{23} proposes that the directional weighted matrix in the directional total variation (DTV)\cite{24} model can be combined with the STV model. The improved model is particularly sensitive to unidirectional images and can effectively remove noise, which is called the direction-guided structure tensor total variation (DSTV)\cite{23} model. But the application scenarios of this model are not very broad, which prompts us to find a proper way to distinguish the gradients in different directions in the STV model and fully explore the underlying image structures.

The STV model utilizes the discrete gradient of the image in the process of calculating the image structure tensor. Therefore, inspired by the  weighted matrix in the ATV model given by Pang et al.\cite{25}, we also consider using the anisotropic matrix to apply different weights to the discrete gradient operator in the $x$-axis and $y$-axis directions in the STV model. As a result, the STV model weighted by an anisotropic diffusion matrix, called weighted structure tensor total variation (WSTV) is given, which penalizes the weighted gradient magnitude before penalizing the  eigenvalues  of the structure tensor. This allows the image gradients in each channel in STV to diffuse along the tangent direction of local features, thereby effectively capturing the local information of the image.

\subsection{Contributions}
In this work, we give a weighted structure tensor total variation model for the image denoising problem. The main contributions of this paper are
\begin{itemize}
	\item[$\bullet$] The WSTV regularization term is proposed for the image denoising minimization model. Compared to the  STV model, a weighted matrix is introduced to achieve varying degrees of punishment for discrete gradients in the $x$ and $y$ directions. The WSTV model  couples the local features of the image and maintains its details well.
	\item[$\bullet$] The WSTV optimization model is solved by an effective first-order algorithm. We show the equivalence of the original model and its dual problem. Also, we give an upper bound of the Lipschitz constant of the dual problem so that the classical FISTA  can be applied to the dual problem efficiently. Theoretically, we prove  the properties of the dual problem,  which guarantee the $O(1/i^2)$ convergence rate of the algorithm.

	\item[$\bullet$] Numerical experiments demonstrate that the proposed WSTV method is superior when compared to a class of  TV-based methods. For most of the grayscale and colored image denoising problems in Section \ref{Numerical experiment}, the values of PSNR and SSIM produced by WSTV are the best.

\end{itemize}

The remainder of this paper is organized as follows. Some elementary knowledge related to image denoising and several denoising models are briefly introduced in Sec.~\ref{Preliminaries}. In Sec.~\ref{sect3}, the WSTV model is given. In addition, we show the computation method that applies the fast gradient projection algorithm to the dual problem of the WSTV optimization model. In Sec.~\ref{Numerical experiment}, we display and discuss the results of numerical experiments on images with different Gaussian noises. Finally,  we conclude the paper in Sec.~\ref{Conclusion}.
\section{Preliminaries}\label{Preliminaries}
In this section, we introduce the symbols and concepts used thereafter. The general linear image restoration models as well as TV and STV-based denoising models are reviewed.

\subsection{Notations}

Unless otherwise specified, vectors and matrices are represented by using lowercase and uppercase bold letters, respectively. We consider the vector-valued image $\boldsymbol{u}=(\boldsymbol{u}_{1}, \boldsymbol{u}_{2}, \cdots, \boldsymbol{u}_{M})$ of ${M}$ channels for convenience, using ${N}$ to represent the total number of pixels in each channel and stacking the channels to obtain  $\boldsymbol{u} \in \mathbb{R}^{{N} {M}}$.
Define $\Upsilon\in \mathcal{Y}\triangleq\mathbb{R}^{M \times N\times 2}$, then the $\ell_{2,1}$ norm  of $\Upsilon$  is
$$\|\Upsilon \|_{2,1}=\sum_{i=1}^{M} \sum_{j=1}^N \sqrt{\sum_{s=1}^{2} \Upsilon_{i, j,s}^{2}}.$$
Apart from this $\ell_{2,1}$ norm, for  $\Lambda =[\Lambda_{1}^{T}, \ldots, \Lambda_{M}^{T}]^{T} \in \mathcal{Y}$,  the norm $\|\Lambda\|_{q, p}$ for any $q\neq2$ or $p\neq1$ represents
$$\|\Lambda\|_{q, p}=(\sum_{i=1}^{M}\|\Lambda_{i}\|_{\mathcal{S}_{p}}^{q})^{1 / q},$$
where $\| \cdot \|_{\mathcal{S}_{p}}$ represents the $p$-order Schatten norm.
The norm $\|\cdot\|_{q, p}$ is also denoted as $\ell_{q}\mbox{-}\mathcal{S}_{p}$ \cite{22}.

Let $\Upsilon, \Lambda \in \mathcal{X}\triangleq\mathbb{R}^{N \times LM\times 2}$, with  $\Upsilon_{n}, \Lambda_{n} \in \mathbb{R}^{LM \times 2}$ $(\forall n=1,2, \ldots, N)$. Then the inner product $\langle\cdot, \cdot\rangle_\mathcal{X}$ and norm $\|\cdot\|_{\mathcal{X}}$ can be denoted by
$$\langle\Lambda , \Upsilon\rangle_{\mathcal{X}}=\sum_{n=1}^{N} \operatorname{trace}(\Upsilon_{n}^{T} \Lambda_{n})$$
and
$$\|\Upsilon\|_{\mathcal{X}}=\sqrt{\langle\Upsilon, \Upsilon\rangle_{\mathcal{X}}}=(\sum_{n=1}^{N}\|\Upsilon_{n}\|_{F}^{2})^{\frac{1}{2}},$$
where $\operatorname{trace}(\cdot)$ is the trace operator of a matrix and $\|\cdot\|_{F}$ is the Frobenius norm. 
The symbol  ${B}_{\infty, q}$  denotes  $\ell_{\infty}\mbox{-}\mathcal{S}_{q}$ unit-norm ball, which is defined as
$${B}_{\infty, q}=\{\boldsymbol{\Phi} \in \mathcal{X}:\|\boldsymbol{\Phi}(n)\|_{\mathcal{S}_{q}} \leq 1, \forall n=1, \ldots, N\},$$ where $\boldsymbol{\Phi}(n)$ represents the $n$-th submatrix of $\boldsymbol{\Phi}$.
\subsection{Image restoration model}
To address the ill-posed problem of image denoising, a basic linear image restoration model is given as
\begin{equation}
\boldsymbol{f}(x)=\boldsymbol{Au}(x)+\boldsymbol{n}(x)\label{a1},
\end{equation}
where $\boldsymbol{u}(x)=(\boldsymbol{u}_{1}(x), \boldsymbol{u}_{2}(x), \cdots, \boldsymbol{u}_{M}(x)): \Omega \rightarrow \mathbb{R}^{{M}}$ is the expected restored image, and $\boldsymbol{f}(x): \Omega \rightarrow \mathbb{R}^{{M}}$ denotes the degraded image. Here $\boldsymbol{A}$ is some linear bounded irreversible operator mapping from one function space to another. For example, the linear operator $\boldsymbol{A}$ is the identity operator for image denoising. The $\boldsymbol{n}(x)$ denotes additional Gaussian noise, and $\Omega\subset\mathbb{R}^{{2}}$ denotes the ${2}$-dimensional image domain.  
The dimension ${M}$ represents the number of image channels.
 In order to recover the unknown image $\boldsymbol{u}$ from Eq.~(\ref{a1}), a general model for such an inverse problem takes the following form
\begin{equation}
\underset{\boldsymbol{u}}{\min } \frac{1}{2}\|\boldsymbol{Au}-\boldsymbol{f}\|_{2}^{2}+\tau R(\boldsymbol{u}). \label{a2}
\end{equation}
The former term of the model is the data fidelity term used to maintain the image structures, while the latter is a regularization term. The regularization parameter $\tau\ge0$ is used to balance both of the above.

Generally, the total variation of noisy images is larger than that of non-noisy images. Therefore, minimizing the total variation can eliminate noise. Image denoising based on total variation can be summarized as the following minimization problem
\begin{equation}
\underset{\boldsymbol{u}}{\min }\frac{\tau}{2}\|\boldsymbol{u}-\boldsymbol{f}\|_{2}^{2}+\|\nabla \boldsymbol{u}\|_{2,1},\label{a3}
\end{equation}
where $\nabla \boldsymbol{u}$ represents the discrete gradient obtained through the forward difference operator\cite{26}. The $\ell_{2,1}$ norm of $\nabla \boldsymbol{u}$ can be calculated by
$\|\nabla \boldsymbol{u} \|_{2,1}=\sum_{i=1}^{NM} \sqrt{\sum_{s=1}^{2} \nabla \boldsymbol{u}_{i, s}^{2}}$ since $\boldsymbol{u}$ has been defined as a vector.

\subsection{Anisotropic total variation (ATV)}
Actually, most TV-based models that we investigate are isotropic. The same weight is given when we calculate the discrete gradient of the model, which results in the same penalty  applied to both the horizontal and vertical directions. However, this setting cannot  well address the local features in the process of restoring images. Consequently, Pang et al. proposed a new anisotropic total variation (ATV)\cite{25}. In order to keep consistent with the formulae in Ref.~\citenum{25}, we define the grayscale image $\boldsymbol{u}\in\mathbb{R}^{M\times N}$, i.e., $\boldsymbol{u}$ is a matrix here.
We clarify that the symbol $\boldsymbol{u}$ refers to the vector-valued image beyond this subsection. 
The ATV model can diffuse along the tangent direction of local features in the image by using an anisotropic weighted matrix, which assigns different weights to discrete gradients in the $x$ and $y$ directions. Specifically, it can  be constructed by
\begin{equation}
\underset{\boldsymbol{u}}{\min } \frac{\tau}{2}\|\boldsymbol{u}-\boldsymbol{f}\|_{2}^{2}+\|\mathcal{W}\nabla \boldsymbol{u}\|_{2,1},\label{a4}
\end{equation}
where  $\nabla \boldsymbol{u}\in \mathcal{Y}$.
The anisotropic weighted matrix $\mathcal{W}$ is defined as
\begin{equation}
\mathcal{W}=\left[\begin{array}{cc}
\boldsymbol{w}_{1} & 0 \\
0 & \boldsymbol{w}_{2}
\end{array}\right]=\left[\begin{array}{cc}
\frac{1}{1+\kappa\left|\boldsymbol{G}_{\hat{\sigma}} * \nabla_{x} \boldsymbol{f}\right|} & 0 \\
0 & \frac{1}{1+\kappa\left|\boldsymbol{G}_{\hat{\sigma}} * \nabla_{y} \boldsymbol{f}\right|}
\end{array}\right],
\label{a5}
\end{equation}
where $\kappa$ and $ \hat\sigma $ are parameters. The symbol $\boldsymbol{G}_{\hat{\sigma }}$ represents the Gaussian convolution of variance $\hat\sigma$ to reduce the influence of noise. The model degenerates into TV when the operators $\boldsymbol{w}_1$ and $\boldsymbol{w}_2$ take the same constant. Since the model (\ref{a4}) is a non-smooth convex optimization problem, the alternating direction method of multipliers (ADMM)  can be used to solve it.

\subsection{Structure tensor total variation (STV)}

The two TV-based models mentioned above penalize the  gradient magnitude, which is completely localized. To address this limitation, Lefkimatis et al. proposed structure tensor total variation (STV)\cite{22}, where the information available in the local neighborhood of each point in the image domain is taken into account. Then the   structure tensor square root of the eigenvalues of the image is penalized. Therefore, the obtained regularizer exhibits semi-local behavior. 

Any pixel  ${i}$  of a vector-valued image  $\boldsymbol{u}$  has a Jacobian matrix, which is defined as
\begin{equation}
( \boldsymbol{J} \boldsymbol{u})(i)=[\nabla \boldsymbol{u}_{{1}}(i), \nabla \boldsymbol{u}_{{2}}(i),\ldots ,\nabla \boldsymbol{u}_{{M}}(i)]^{{T}},\label{a6}
\end{equation}
where $\nabla \boldsymbol{u}_{1}(i)$ represents discrete gradient\cite{26}. Define the structure tensor of image $\boldsymbol{u}$ at pixel $i$ as
\begin{equation}
({S}_{K} \boldsymbol{u})(i)= {K} *[(\boldsymbol{Ju})(i)^{{T}}(\boldsymbol{Ju})(i)],\label{a7}
\end{equation}
where $ {K}$ is a Gaussian convolution kernel, and $*$ represents the  convolution operation. Define $\lambda^{+}=\lambda^{+}(S_{K} \boldsymbol{u}(i))$ and $\lambda^{-}=\lambda^{-}(S_{K} \boldsymbol{u}(i))$ as the maximum and minimum eigenvalues of the structure tensor $S_{K}\boldsymbol{u}(i)$ at any point $i$ in image $\boldsymbol{u}$. From Ref.~\citenum{22}, we know that the $\ell_{p}$ norms of the square root of the eigenvalues of the structure tensor are more effective in measuring local image variation. Naturally, the STV regularizer can be defined as
\begin{equation}
\operatorname{STV}_{{p}}(\boldsymbol{u})=\sum_{{i}=1}^{{N}}\|(\sqrt{\lambda_{{i}}^{+}}, \sqrt{\lambda_{{i}}^{-}})\|_{{p}},{p}\ge1 \label{a8}.
\end{equation}
	
Because of the nonlinearity of the operator and the existence of the convolution kernel $K$ in $S_{K}\boldsymbol{u}(i)$, Lefkimmiatis et al.\cite{22} proposed an alternative formulation of $S_{K}\boldsymbol{u}(i)$ and named it  patch-based Jacobian $\boldsymbol{J}_{K}: \mathbb{R}^{{N} {M}} \rightarrow  \mathbb{R}^{{N} \times{LM} \times 2}$. For any pixel $i$ in the image $\boldsymbol{u}$, we have
\begin{equation}
(\boldsymbol{J}_{K}\boldsymbol{u})(i)=[\tilde\nabla \boldsymbol{u}_{1}(i), \tilde\nabla \boldsymbol{u}_{2}(i), \ldots  ,\tilde\nabla \boldsymbol{u}_{{M}}(i)
]^{{T}},\label{a9}
\end{equation}
where $\tilde\nabla \boldsymbol{u}_{{1}}(i)=[({\Psi}_{1} \nabla \boldsymbol{u}_{{1}})(i),({\Psi}_{2} \nabla \boldsymbol{u}_{{1}})(i), \ldots,({\Psi}_{{L}}\nabla \boldsymbol{u}_{{1}})(i)]$.  In addition, ${L}=(2 {L}_{{{K}}}+1)^{2}$ indicates the number of all elements in the convolution kernel $K$. The weighted translation operator ${\Psi}_{l}(l=1, \ldots, L)$ is defined as $ ({\Psi}_{l}\nabla \boldsymbol{u}_{{m}})(i)=\sqrt{K[g_{l}]}(\nabla\boldsymbol{u}_{{m}})(x_{i}-g_{l})$, where $g_{l}\in\{-L_{K}, \ldots, L_{K}\}^{2}$  is the shift amount. Based on the above definition, the discrete structure tensor of $\boldsymbol{u}$ evaluated at the pixel location $i$ in terms of the patch-based Jacobian can be written as
\begin{equation}
(S_{K}\boldsymbol{u})(i)=((\boldsymbol{J}_{K} \boldsymbol{u})(i))^{T}(\boldsymbol{J}_{K}\boldsymbol{u})(i).\label{a10}
\end{equation}
Lefkimmiatis et al. \cite{22} have also proved that the singular values of $(\boldsymbol{J}_{K} \boldsymbol{u})(i)$ are equivalent to the square root of the eigenvalues of $(S_{K} \boldsymbol{u})(i)$. Then the STV regularizer in Eq.~$(\ref{a8})$ can be redefined as
\begin{equation}
\operatorname{STV}_{p}(\boldsymbol{u})=\sum_{i=1}^{N}\|(\boldsymbol{J}_{K} \boldsymbol{u})(i)\|_{\mathcal{S}_{p}}=\|{{\boldsymbol{J}_{K}} \boldsymbol{u}}\|_{1, p}.\label{a11}
\end{equation}
\section{The proposed model and algorithm}\label{sect3}	
\subsection{The weighted structure tensor total variation (WSTV) model}
		
In order to fully explore the local features of the image, we employ the anisotropic weighted matrix  proposed in  the ATV model  to refine the patch-based Jacobian operator of the STV model. The new weighted patch-based Jacobian $	(\widehat{\boldsymbol{J}_{K}} \boldsymbol{u})(i)$ can be represented as
\begin{equation}
(\widehat{\boldsymbol{J}_{K}} \boldsymbol{u})(i)=[
\widehat\nabla  \boldsymbol{u}_{1}(i),\widehat\nabla  \boldsymbol{u}_{2}(i),  \ldots ,\widehat\nabla  \boldsymbol{u}_{M}(i)]^{T},\label{a12}
\end{equation}
where
$\widehat\nabla  \boldsymbol{u}_{1}(i)=[({\Psi}_{1}\mathcal{W} \nabla \boldsymbol{u}_{{1}})(i),({\Psi}_{2}\mathcal{W} \nabla \boldsymbol{u}_{{1}})(i),\ldots,({\Psi}_{L}\mathcal{W} \nabla \boldsymbol{u}_{{1}})(i)]$.  
Here $\mathcal{W}$ represents the anisotropic weighted matrix defined in  Eq.~(\ref{a5}). Furthermore, the existence of $\mathcal{W}$ will not destroy the convexity of the regularizer. The WSTV regularizer can be defined as
\begin{equation}
\operatorname{WSTV}_{p}(\boldsymbol{u})=\sum_{i=1}^{N}\|(\widehat{\boldsymbol{J}_{K}}\boldsymbol{u})(i)\|_{\mathcal{S}_{p}}=\|\widehat{\boldsymbol{J}_{K}}\boldsymbol{u}\|_{1, p}.\label{a13}
\end{equation}
According to the objective function given in Eq.~(\ref{a2}), the problem to be solved is
\begin{equation}
\underset{\boldsymbol{u} \in \mathcal{C}}{\min } \frac{1}{2}\|\boldsymbol{u}-\boldsymbol{f}\|_{2}^{2}+\tau\|\widehat{\boldsymbol{J}_{K}} \boldsymbol{u}\|_{1, p}.\label{a14}
\end{equation}	
Since the WSTV regularizer is convex, the minimizer $\hat{\boldsymbol{u}}$ is the proximal point generated by the proximal operator associated with the regularizer $\operatorname{WSTV}_{p}$ at $\boldsymbol{f}$, i.e.,
\begin{equation}
\hat{\boldsymbol{u}}=\operatorname{prox}_{\tau\operatorname{WSTV}_{p}(\boldsymbol{u})}(\boldsymbol{f}):=\underset{\boldsymbol{u} \in \mathcal{C}}{\arg \min } \frac{1}{2}\|\boldsymbol{u}-\boldsymbol{f}\|_{2}^{2}+\tau\|\widehat{\boldsymbol{J}_{K}} \boldsymbol{u}\|_{1, p},\label{a15}
\end{equation}
where  $\mathcal{C}$  is a convex set that represents additional constraints on the solution. In this paper, we consider the box constraints $(\mathcal{C}=\mathbb{R}^{NM}$  in unconstrained case).
			
\subsection{Algorithm for WSTV}
Since the suggested WSTV regularizer is not differentiable, we consider a dual formulation for the optimization problem (\ref{a15}). We show the computation framework in two steps. First, the dual model as well as the equivalence between the primal and dual problems are deduced. The convexity, Lipschitz continuity, and solvability of the dual model are also proved. Then, we employ the fast gradient projection method to find the optimal solution of the dual model. Moreover, details and the complexity result $O(1/i^2)$ of the algorithm are given.
			
\vspace{2ex}\noindent{\footnotesize\textbf{Lemma 1.}
Let  $p \geq 1$, and let $ q$  be the conjugate exponent of  $p$, i.e., $ \frac{1}{p}+\frac{1}{q}=1$. Then, the mixed vector-matrix norm $ \|\cdot\|_{\infty, q}$  is dual to the mixed vector-matrix norm  $\|\cdot\|_{1, p}${\cite{27}}}.
			
Through Lemma 1  and the property that the dual of the dual norm is the original norm, we are able to redefine the WSTV regularizer in terms of
\begin{equation}
\|{\widehat{\boldsymbol{J}_{K}} \boldsymbol{u}}\|_{1, p}=
\underset{\boldsymbol{\Phi} \in {B}_{\infty, q}}{\max }
\langle\mathbf{\Phi }, \widehat{\boldsymbol{J}_{K}} \boldsymbol{u}\rangle_\mathcal{X},\label{a16}
\end{equation}
where $\boldsymbol{\Phi}$ denotes the variable in the target space $\mathcal{X}\triangleq\mathbb{R}^{N \times LM \times 2}$. 
Then, we can rewrite Eq.~(\ref{a15})  as follows
\begin{equation}
\begin{aligned}
\hat{\boldsymbol{u}}&=\underset{\boldsymbol{u} \in \mathcal{C}}{\arg \min } \frac{1}{2}\|\boldsymbol{u}-\boldsymbol{f}\|_{2}^{2}+\tau \underset{\boldsymbol{\Phi} \in {B}_{\infty, q}}{\max }\langle\mathbf{\Phi }, \widehat{\boldsymbol{J}_{K}} \boldsymbol{u}\rangle_\mathcal{X}\\&=\underset{\boldsymbol{u} \in \mathcal{C}}{\arg \min } \frac{1}{2}\|\boldsymbol{u}-\boldsymbol{f}\|_{2}^{2}+\tau \underset{\boldsymbol{\Phi} \in {B}_{\infty, q}}{\max } \langle\widehat{\boldsymbol{J}_{K}}^{*} \mathbf{\Phi}, \boldsymbol{u}\rangle_{2}.
\end{aligned}\label{a17}
\end{equation}
The adjoint operator $\widehat{\boldsymbol{J}_{K}}^{*}:\mathbb{R}^{N \times L M \times 2} \rightarrow  \mathbb{R}^{NM}$ of the weighted patch-based Jacobian $\widehat{\boldsymbol{J}_{K}}$ is defined as
\begin{equation}
(\widehat{\boldsymbol{J}_{K}}^{*}\boldsymbol{X})(t)=\sum_{l=1}^{L}-\operatorname{\mathbf{div_{(\delta )}}}({\Psi}_{l}^{*}  \boldsymbol{X}(i,h)),\label{a18}
\end{equation}
where $h=(m-1)L+l$ and $t=(m-1)N+n$  with  $1 \leq n \leq N$  and  $1 \leq m \leq M$. The symbol ${\Psi}_{l}^{*}$ is the adjoint of the weighted translation operator ${\Psi}_{l}$, which scans the $\boldsymbol{X}(i,h)$ in column-wise manner. The two-element vector $\boldsymbol{X}(i,h)$ is the $h$th row of the $i$th submatrix of an arbitrary $\boldsymbol{X}\in\mathbb{R}^{N \times LM \times 2}$ [\citenum{22}, Proposition 3.2]. The weighted discrete divergence operator $\operatorname{\mathbf{div_{(\delta )}}}$ is denoted by
\begin{equation}
\operatorname{\mathbf{div_{(\delta )}}} (\cdot)\triangleq \operatorname{\mathbf{div}(\mathcal{W}^{{T}}(\cdot))}= \operatorname{\mathbf{div}(\mathcal{W}(\cdot))}.\label{a19}
\end{equation}
Here $\mathbf{div}$ represents the discrete divergence, which is defined using backward differences\cite{26}.

The Eq.~(\ref{a17})  is rewritten as a minimax problem
\begin{equation}
\underset{\boldsymbol{u} \in \mathcal{C}}{\min }
\underset{\boldsymbol{\Phi} \in {B}_{\infty, q}}{\max } H(\boldsymbol{u}, \boldsymbol{\Phi})\label{a20},
\end{equation}	
where $H(\boldsymbol{u}, \boldsymbol{\Phi})=\frac{1}{2}\|\boldsymbol{u}-\boldsymbol{f}\|_{2}^{2}+\tau\langle\widehat{\boldsymbol{J}_{K}}^{*} \mathbf{\Phi}, \boldsymbol{u}\rangle_{2}$.
Since the function $H(\boldsymbol{u},\boldsymbol{\Phi})$ is  convex with respect to $\boldsymbol{u}$ and concave with respect to $\boldsymbol{\Phi}$, there is a common saddle point $(\hat{\boldsymbol{u}}, \hat{\mathbf{\Phi }})$ that does not change when the minimum and the maximum are interchanged. Specifically, we have
\begin{equation}
\underset{\boldsymbol{u} \in \mathcal{C}}{\min }
\underset{\boldsymbol{\Phi} \in {B}_{\infty, q}}{\max }H(\boldsymbol{u}, \boldsymbol{\Phi })=H(\hat{\boldsymbol{u}}, \hat{\mathbf{\Phi }})=\underset{\boldsymbol{\Phi} \in {B}_{\infty, q}}{\max }
\underset{\boldsymbol{u} \in \mathcal{C}}{\min }\;
H(\boldsymbol{u}, \boldsymbol{\Phi })\label{a21}.
\end{equation}
Maximizing the dual objective function  $d(\boldsymbol{\Phi})=\underset{\boldsymbol{u} \in \mathcal{C}} {\min}\; H(\boldsymbol{u}, \boldsymbol{\Phi})$  in Eq.~(\ref{a21}) is equivalent to minimizing the primal objective function $p(\boldsymbol{u})=\underset{\boldsymbol{\Phi} \in {B}_{\infty, q}}{\max } H(\boldsymbol{u}, \boldsymbol{\Phi})$. As a result, we can  find the minimizer  $\hat{\boldsymbol{u}} $ of  $p(\boldsymbol{u})$ by seeking the maximizer  $\hat{\boldsymbol{\Phi}} $ of  $d(\boldsymbol{\Phi})$.
The minimization problem of ${\boldsymbol{u}}$ can be obtained by expanding the function $ H(\boldsymbol{u}, \boldsymbol{\Phi})$, i.e.,
\begin{equation}
\hat{\boldsymbol{u}}=\underset{\boldsymbol{u} \in \mathcal{C}}{\operatorname{argmin}}\|\boldsymbol{u}-(\boldsymbol{f}-\tau \widehat{\boldsymbol{J}_{K}}^{*} \boldsymbol{\Phi})\|_{2}^{2}-D,\label{a22}
\end{equation}
where $ D $ represents a constant term and the solution of Eq.~(\ref{a22}) is $ \hat{\boldsymbol{u}}=P_{\mathcal{C}}(\boldsymbol{f}-\tau \widehat{\boldsymbol{J}_{K}}^{*} \boldsymbol{\Phi})$. The operator 
$ P_{\mathcal{C}}$  is the orthogonal projection operator on the convex set $\mathcal{C} $. By substituting  $\hat{\boldsymbol{u}}$  into the function  $H(\boldsymbol{u}, \boldsymbol{\Phi})$,  we have
\begin{equation}
d(\boldsymbol{\Phi})=H(\hat{\boldsymbol{u}}, \boldsymbol{\Phi})=\frac{1}{2}\|\boldsymbol{s}-P_{\mathcal{C}}(\boldsymbol{s})\|_{2}^{2}+\frac{1}{2}\|\boldsymbol{f}\|_{2}^{2}-\frac{1}{2}\|\boldsymbol{s}\|_{2}^{2},\label{a23}
\end{equation}
where $\boldsymbol{s}=\boldsymbol{f}-\tau \widehat{\boldsymbol{J}_{K}}^{*} \boldsymbol{\Phi} $. Therefore, the dual problem of the WSTV model is finally presented as  			
\begin{equation}
\hat{\boldsymbol{\Phi}}=\underset{\boldsymbol{\Phi} \in {B}_{\infty, q}}{\arg \max} \ d(\boldsymbol{\Phi})
 =\underset{\boldsymbol{\Phi} \in {B}_{\infty, q}}{\arg \max } \frac{1}{2}\|\boldsymbol{s}-P_{\mathcal{C}}(\boldsymbol{s})\|_{2}^{2}+\frac{1}{2}\|\boldsymbol{f}\|_{2}^{2}-\frac{1}{2}\|\boldsymbol{s}\|_{2}^{2}.\label{a24}
\end{equation}

In the next theorem, we show the properties of the dual problem.

\vspace{2ex}\noindent{\textbf{Theorem 1.}
The dual optimization model  in (\ref{a24}) satisfies the following properties
\begin{itemize}
\item [(1)] The objective function $d(\boldsymbol{\Phi})$  is convex;
\item [(2)] The gradient of $d(\boldsymbol{\Phi})$ is  Lipschitz continuous with Lipschitz  constant  ${L}(d) = 8 \sqrt{2} \tau^{2}$;
\item [(3)] The dual optimization problem is solvable, i.e., $\hat{\boldsymbol{\Phi}}$ is not an empty set.
\end{itemize}}
\vspace{2ex}\noindent{\textbf{Proof.}
(1) Because $H(\boldsymbol{u}, \boldsymbol{\Phi})$ is concave with respect to $\boldsymbol{\Phi}$, then  $\underset{\boldsymbol{u} \in \mathcal{C}}{\sup}\;\{-H(\boldsymbol{u}, \boldsymbol{\Phi})\}$ is a convex function with respect to $\boldsymbol{\Phi}$. That is to say, $d(\boldsymbol{\Phi})=\underset{\boldsymbol{u} \in \mathcal{C}} {\min}\;H(\boldsymbol{u}, \boldsymbol{\Phi})$ is a convex function. The first part of the theorem holds.		
	
(2) Similar to Proposition 4.2 in Ref.~\citenum{22}, we have
\begin{equation}
\begin{aligned}
\|\nabla d(\boldsymbol{a})-\nabla d(\boldsymbol{b})\| &=\tau\|\widehat{\boldsymbol{J}_{K}}P_{\mathcal{C}}(\boldsymbol{f}-\tau {\widehat{\boldsymbol{J}_{K}}}^{*} \boldsymbol{a})-\widehat{\boldsymbol{J}_{K}}P_{\mathcal{C}}(\boldsymbol{f}-\tau {\widehat{\boldsymbol{J}_{K}}}^{*} \boldsymbol{b})\| \\	
& \leq \tau\|{\widehat{\boldsymbol{J}_{K}}}\|\|P_{\mathcal{C}}(\boldsymbol{f}-\tau {\widehat{\boldsymbol{J}_{K}}}^{*} \boldsymbol{a})-P_{\mathcal{C}}(\boldsymbol{f}-\tau{\widehat{\boldsymbol{J}_{K}}}^{*} \boldsymbol{b})\| \\
& \leq \tau\|\widehat{\boldsymbol{J}_{K}}\|\|\tau {\widehat{\boldsymbol{J}_{K}}}^{*}(\boldsymbol{a}-\boldsymbol{b})\| \\		
& \leq \tau^{2}\|\widehat{\boldsymbol{J}_{K}}\|^{2}\|\boldsymbol{a}-\boldsymbol{b}\| \\	
& \leq \tau^{2}\|\mathcal{W}\|^{2}\|\nabla\|^{2}\|T\|^{2}\|\boldsymbol{a}-\boldsymbol{b}\|,\nonumber
\end{aligned}
\end{equation}
where  $T=\sum_{l=1}^{L}({\Psi}_{l}^{*} {\Psi}_{l})$. It is proved in Ref.~\citenum{29}  that $ \|\nabla\|^{2} \leq 8 $ and  $\|T\|^{2}\leq\sqrt{2}$, and further we have $ \|\mathcal{W}\|^{2} \leq \max \{\boldsymbol{w}_{1}^{2}, \boldsymbol{w}_{2}^{2}\}^{2} \leq 1 $. These indicate that  ${L}(d) = 8 \sqrt{2} \tau^{2}$.

(3) Since the objective function $d(\boldsymbol{\Phi})$ is  smooth and the  constraint $\boldsymbol{\Phi}\in{B}_{\infty, q}$ is compact, the dual problem is solvable. The proof is completed.}

Because the dual problem is smooth and convex,  we are able to apply the fast gradient projection method (FGP)\cite{29} to solve it. In each iteration, a transitional point is obtained by the gradient projection method. Then, a new iteration point is generated by the accelerated technique in FISTA from the transitional point and the previous point. The detailed computational process involves the gradient of the objective function in (\ref{a24}), the step size of the first-order algorithm, and the projection onto the constraint ${B}_{\infty, q}$.

The gradient of the objective function in the dual problem is
$$\nabla d(\boldsymbol{\Phi })=\tau \widehat{\boldsymbol{J}_{K}}P_{\mathcal{C}}(\boldsymbol{f}-\tau \widehat{\boldsymbol{J}_{K}}^{*} \boldsymbol{\Phi}).$$
For the step size along the gradient ascent direction, a constant of $1/L(d)$ is chosen to ensure the convergence of the FGP algorithm. Here, $L(d)$ represents the Lipschitz constant whose upper bound  $\leq 8 \sqrt{2} \tau^{2}$ is proved in Theorem 1.
In terms of the projection of a matrix on the unit normal ball ${B}_{\infty, q}$, the authors of Ref.~\citenum{22} proved that it can be calculated with the aid of the SVD decomposition of the matrix. Below we give the explicit steps for calculating the projection of  $\boldsymbol{\Phi}(n)$ on ${B}_{\infty, q}$.
\begin{enumerate}
  \item [(1)] Conduct SVD to  $\boldsymbol{\Phi}(n)$, i.e., $\boldsymbol{\Phi}(n)=\boldsymbol{U}\boldsymbol{\Sigma} \boldsymbol{V}^{T}$, where $\boldsymbol{\Sigma}=\operatorname{diag}(\sigma_{1}, {\sigma}_{2})$.
  \item [(2)] Calculate the projection of $\boldsymbol{\Phi}(n)$ onto the unit-norm ball  ${B}_{\mathcal{S}_{q}}$, i.e., $P_{{B}_{\mathcal{S}_{q}}}(\boldsymbol{\Phi}(n))=\boldsymbol{U} \boldsymbol{\Sigma}_{q} \boldsymbol{V}^{T}$, where $\boldsymbol{\Sigma}_{q}=\operatorname{diag}(\boldsymbol\sigma_{q})$ and $\boldsymbol\sigma_{q}$ are the projected singular values of $\boldsymbol{\Sigma} $ onto the $\ell_{q}$ unit-norm ball ${B}_{q}=\{\boldsymbol{\sigma} \in \mathbb{R}_{+}^{N}:\|\boldsymbol\sigma\|_{q} \leq 1\}$.
  \item [(3)]  Compute $\boldsymbol{V}$ and $\boldsymbol{\Sigma}$ using the eigenvalue decomposition of $(\boldsymbol{\Phi}(n))^{T}(\boldsymbol{\Phi}(n))$.
 \item [(4)] The projection can be constructed by $ P_{{B}_{\mathcal{S}_{q}}}(\boldsymbol{\Phi}(n))=\boldsymbol{\Phi}(n) \boldsymbol{V} \boldsymbol{\Sigma}^{+} \boldsymbol{\Sigma}_{q} \boldsymbol{V}^{T}$. $\boldsymbol{\Sigma}^{+}$  is the pseudoinverse matrix of  $\boldsymbol{\Sigma}$.
\end{enumerate}
We only consider the case of $ q=\infty$ in this work. Therefore, we have
\begin{equation}
P_{{B}_{\mathcal{S}_{\infty}}}(\boldsymbol{\Phi}(n))=\boldsymbol{\Phi}(n) \boldsymbol{V} \boldsymbol{\Sigma}^{+}\operatorname{diag}(\min (\boldsymbol{\sigma}(\boldsymbol{\Phi}(n)), \mathbf{1})) \mathbf{V}^{T},\label{a25}
\end{equation}
where $\mathbf{1}$ is a vector with all elements set to one and $\boldsymbol\sigma(\boldsymbol{\Phi}(n)) \in \mathbb{R}_{+}^{2}$. Discussions about $q=1,2$ and other situations can be found in Refs.~\citenum{27} and ~\citenum{31}. Finally, we present the complete computational process in Algorithm \ref{alg1}.
\begin{algorithm}
\caption{: FGP algorithm for WSTV-based denoising.}
\label{alg1}
\begin{algorithmic}[1]
\REQUIRE  $\boldsymbol{f}$, $\tau>0$, $p=1$, $P_{\mathcal{C}}$.
\STATE Initialize: $\boldsymbol{\Phi}_{0}=\boldsymbol0\in \mathbb{R}^{N \times LM \times 2}$, $t_{1}=1$, $i=1$.
\WHILE  {stopping  criterion is not satisfied}
\STATE $z=P_{\mathcal{C}}(\boldsymbol{f}-\tau \widehat{\boldsymbol{J}_{K}}^{*} \boldsymbol{\Phi}_{i-1})$
\STATE $\boldsymbol{\Phi}_{i}=P_{{B}_{\infty, \infty}}(\boldsymbol{\Phi}_{i-1}+\frac{1}{8 \sqrt{2} \tau} \widehat{\boldsymbol{J}_{K}}z)$
\STATE $t_{i+1} = \frac{1+\sqrt{1+4t_{i}^{2} }}{2}$
\STATE $\boldsymbol{\Phi}_{i+1} = \boldsymbol{\Phi}_{k}+(\frac{t_{i}-1}{t_{i+1}})(\boldsymbol{\Phi}_{i}-\boldsymbol{\Phi}_{i-1}) $
\STATE $i=i+1$
\ENDWHILE
\ENSURE $\boldsymbol{\hat{u}}=P_{C}(\boldsymbol{f}-\tau {\widehat{\boldsymbol{J}_{K}}}^{*}\boldsymbol{\Phi}_{i})$.
\end{algorithmic}
\end{algorithm}

The above FGP method can be regarded as a special case of the FISTA approach for constrained optimization problems. Thus, the properties of the dual problem proved in Theorem 1 guarantee that [\citenum{30}, Theorem 4.4]
\begin{equation}
d(\boldsymbol{\Phi}_{i})-d(\boldsymbol{\Phi}^{*}) \leq \frac{2 L(d)\|\boldsymbol{\Phi}_{0}-\boldsymbol{\Phi}^{*}\|^{2}}{(i+1)^{2}},\label{a26}
\end{equation}
where $\boldsymbol{\Phi}^{*}$ represents the optimal point. Hence, the complexity result for Algorithm \ref{alg1} is $O(1/i^2)$.
					
\section{Numerical experiment}\label{Numerical experiment}
	
\begin{figure}[H]
\centering
\includegraphics[width=1\linewidth]{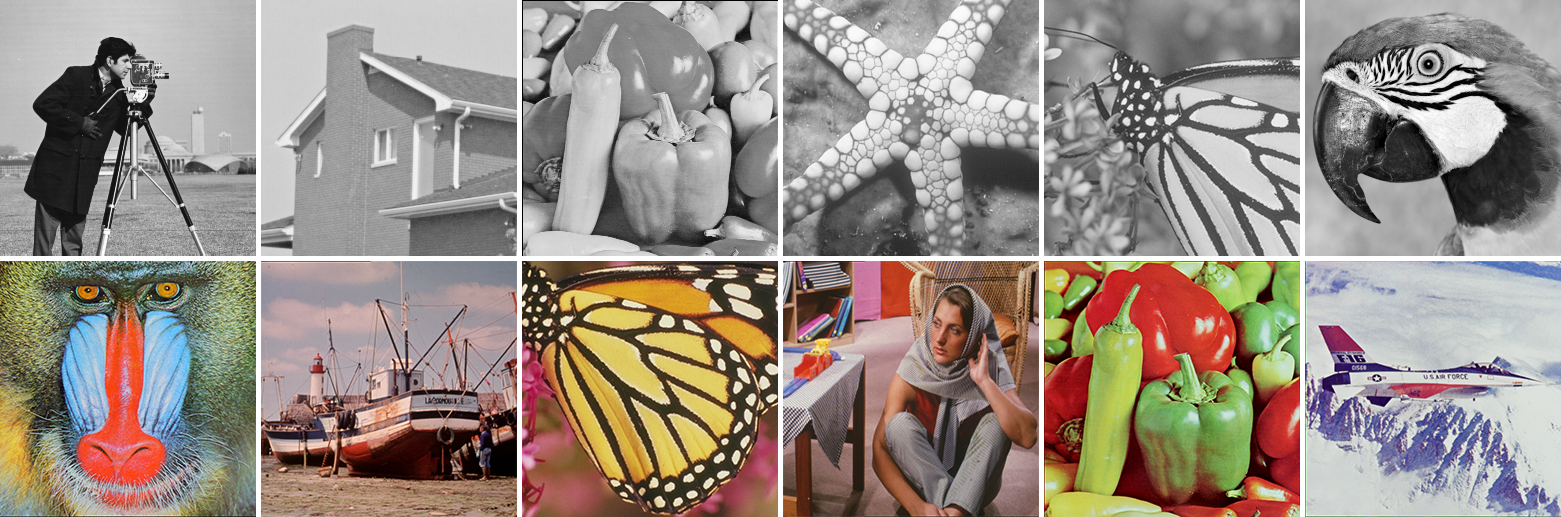}
\caption {Twelve 256 $\times$ 256 test images. From left to right and top to bottom: Cameraman, House, Peppers, Starfish, Butterfly, Parrot, Baboon, Boats, Colored butterfly, Babara, Colored peppers and Airplane.}
\label{fig1}
\end{figure}

In this section, the experimental results are presented to assess the performance of the proposed WSTV model. We mainly focus on the comparisons between the WSTV model and other TV-based models, such as the TV, ATV, VTV, and STV  models. All of the numerical experiments are carried out using MATLAB (R2019a) on a Windows 10 (64-bit) desktop computer  powered by an Intel Core i7 CPU running at 2.70 GHz and 8.0GB of RAM. We employ the structural similarity index (SSIM) and peak signal-to-noise ratio (PSNR) to evaluate the recovered quality of each method. In addition, the intensity of the image involved in the experiment is normalized to the range $[0, 1]$. Furthermore, all numerical methods will be stopped when the relative difference between two successive iterations satisfies
\begin{equation}
\frac{\|u^{i+1}-u^{i}\|_{2}}{\|u^{i}\|_{2}} \leq 10^{-5}\nonumber
\end{equation}
or after reaching a maximum number of iterations. The maximum number of iterations for TV, ATV, and VTV is set to 500, while STV and WSTV have a maximum of 100.
The test images in Figure \ref{fig1} are public domain images of digital image processing.
			
The regularization parameter $\tau$, which controls how much filtering is introduced by the regularization term, has a significant impact on the effectiveness of the model restoration. We adjust the regularization parameters of the five models to empirically make the numerical performance of each method achieve its best.
			
For the tested  grayscale images, Cameraman, House, Peppers, Starfish, Butterfly,  and Parrot. Cameraman  and House  are mainly composed of smooth regions, while Butterfly and Parrot have more edges and textures, respectively. We generate test data by introducing Gaussian noise with standard deviations of $\sigma = \{0.01, 0.05, 0.1, 0.15\}$. Table \ref{tab1} reveals that four models all exhibit strong restoration performance at low noise levels $\sigma = \{0.01, 0.05\}$. As we add very little noise, most details of the image can be preserved during the process of  image restoration. Our model produces higher PSNR and SSIM than those TV-based models at high Gaussian noise levels $\sigma = \{0.1, 0.15\}$.

\begin{figure}[htbp]
\centering
\includegraphics[width=1\linewidth]{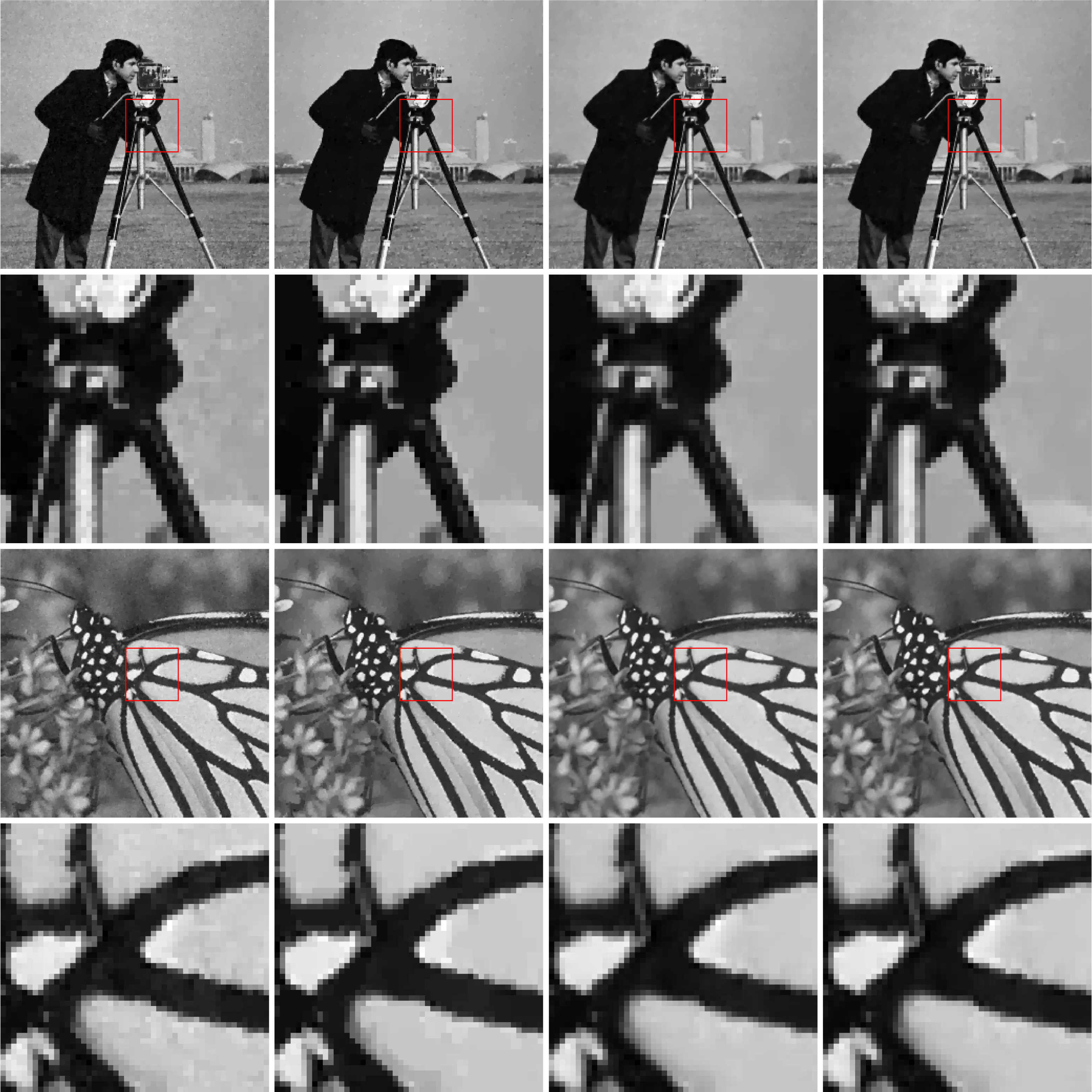}
\caption{\label{fig2}Compare the denoising effects of various models on grayscale images at the Gaussian noise level $\sigma=0.05$. From left to right are the image restoration results obtained using TV, ATV, STV, and the proposed model (WSTV). The first row of each image represents the complete restored image. The second row of each image is a part of the restored image.}
\end{figure}

To effectively demonstrate the performance of the mentioned models for restoring images, we enlarge a portion of the restored image, as shown in Figure \ref{fig2}. We observe that  TV  still retains a large amount of noise in the smooth part of  Cameraman, and the edges are relatively blurry. ATV  has obvious staircase effects but only retains a small amount of noise, while  STV and our WSTV model perform well.

\begin{figure}[htbp]
\centering
\includegraphics[width=1\linewidth]{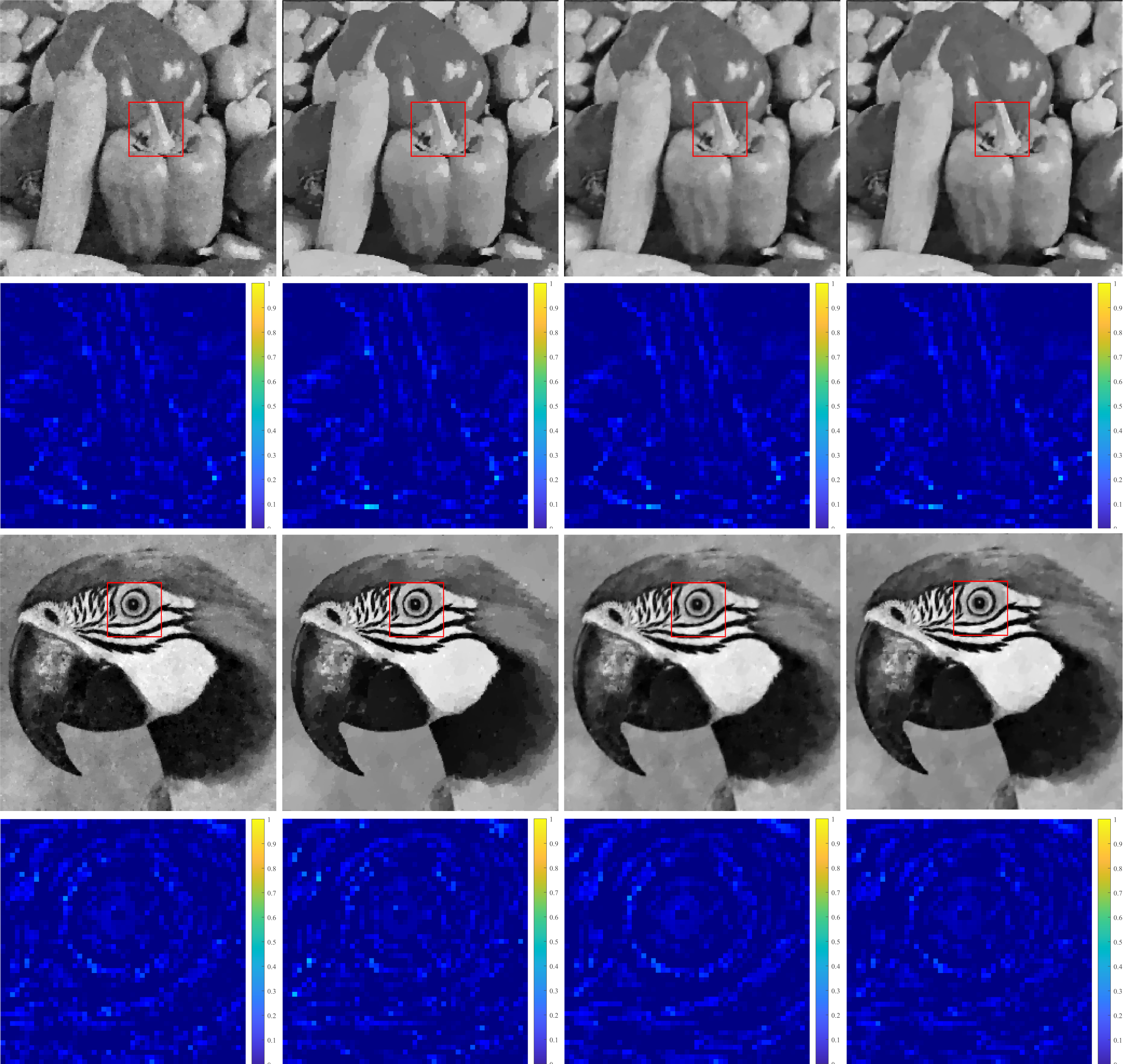}
\caption{Compare the denoising effects of various models on grayscale images at the Gaussian noise level $\sigma=0.1$. From left to right are the image restoration results obtained using TV, ATV, STV, and the proposed model (WSTV). The first row of each image represents the complete restored image. The second row of each image is the difference between the restored image and the original image. The colorbar displays more efficient restorations if the color is more shaded.}
\label{fig3}
\end{figure}

In fact, our model is more effective at restoring details such as edges and corners compared with STV. We display the differences between the restored image and the original image in Figure \ref{fig3}. It is worth noting that the darker the image color, the better the restoration effect. We note the differences exist in Peppers, and there are still many areas with higher brightness in TV, such as the upper right corner of the image. The brightness of the lower edge of ATV is obvious, which indicates that ATV excessively magnifies the local features of the image. Obviously, our model is the darkest, with neat image edges and sufficient details, which means that our method is more robust than other models in restoring the degraded images.
 
For colored image denoising problems, the WSTV method is compared with the vector-extended VTV \cite{32} of TV. Table \ref{tab2} reports the PSNR and SSIM values obtained by three models under different Gaussian noise levels. Similar to gray image denoising, all three models exhibit good restoration performance at lower levels of Gaussian noise. In the images of Barbara and Peppers, STV and WSTV generate the same SSIM, which we believe is caused by the addition of less noise. At high noise levels, our model is significantly superior to the first two.

\begin{figure}[H]
\centering
\includegraphics[width=1\linewidth]{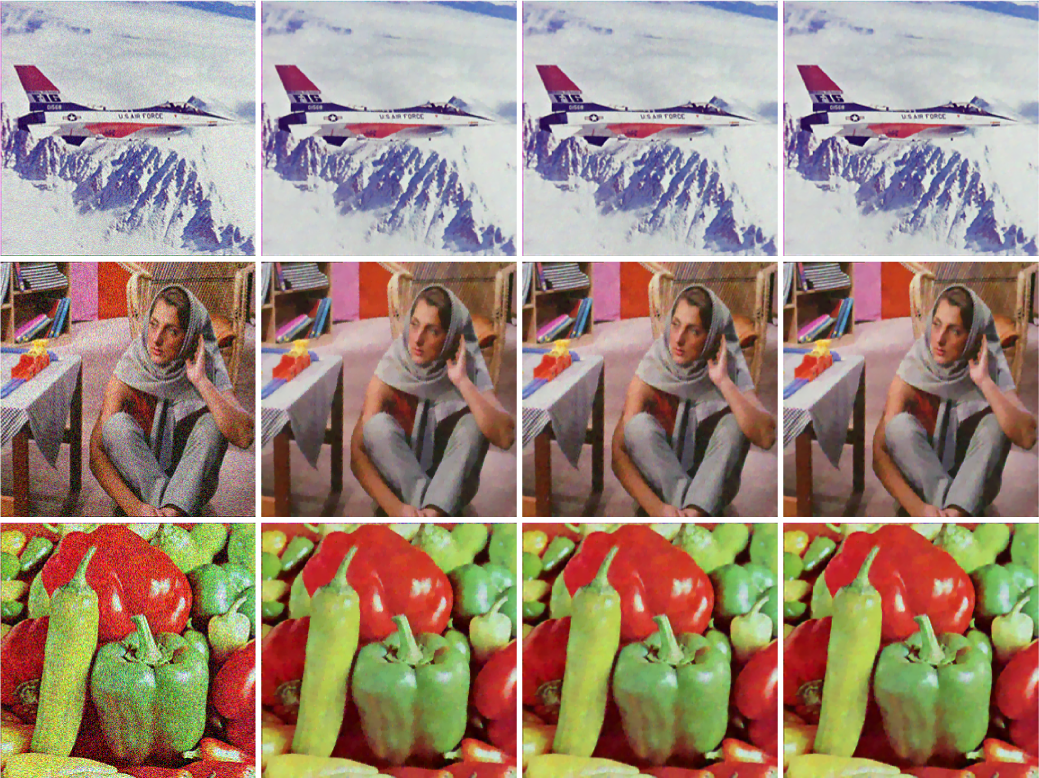}
\caption{\label{fig5}Compare the denoising effects of three models under various noise levels. The Gaussian noise levels from top to bottom are $\sigma=\{0.05, 0.1, 0.15\}$. The restoration results using VTV, STV, and WSTV are shown from the second column on the left to the last column on the right. The first column is the initial noise image.}
\end{figure}

A special case is Baboon, where  STV produces better PSNR values than our model at noise levels of $\sigma=\{0.01, 0.05\}$. In fact, we observe that the details of  Baboon are very abundant, which means that the local information of the image is more complex. Weighting the discrete gradient may have an impact on the information in the pixel domain of images and then reduce the effectiveness of image restoration if the diagonal elements of the anisotropic diffusion matrix are large. Moreover, we show the partial restored images of the three denoising models  with Gaussian noise levels of $\sigma=\{0.05, 0.1, 0.15\}$ in Figure \ref{fig5}. It explicitly demonstrates the superiority of our model in color image denoising for preserving more details of the images.
								
\begin{table}[H]
\centering
\caption{Compare the PSNR and SSIM of different models in grayscale image denoising.}
\resizebox{\textwidth}{!}{
\begin{tabular}{lllllllll}
\toprule\toprule
Noise  & 0.01    &       & 0.05    &         & 0.10    &         &0.15&                       \\ \midrule[1pt]
Image  &  \multicolumn{8}{l}{Cameraman}                                                                                                              \\
Models & PSNR             & SSIM            & PSNR             & SSIM            & PSNR             & SSIM            & PSNR             & SSIM            \\
TV     & 40.6641          & 0.9639          & 30.9808          & 0.8224          & 27.4314          & 0.7269          & 25.7675          & 0.7306          \\
ATV    & 40.5000          & 0.9649          & 31.3241          & 0.8776          & 27.9608          & 0.8076          & 26.1333          & 0.7596          \\
STV    & 41.8871          & 0.9762          & 31.4690          & 0.8893          & 28.0404          & 0.8077          & 26.1356          & 0.7557          \\
WSTV   & \textbf{41.9646} & \textbf{0.9764} & \textbf{31.7752} & \textbf{0.8900} & \textbf{28.4345} & \textbf{0.8211} & \textbf{26.5238} & \textbf{0.7728} \\
Image  &  \multicolumn{8}{l}{Butterfly}                                                                                                              \\
Models & PSNR             & SSIM            & PSNR             & SSIM            & PSNR             & SSIM            & PSNR             & SSIM            \\
TV     & 39.3814          & 0.9757          & 30.7683          & 0.8773          & 27.2040          & 0.8035          & 25.5354          & 0.7922          \\
ATV    & 40.8081          & 0.9774          & 31.3426          & 0.9098          & 27.5779          & 0.8450          & 25.4152          & 0.7886          \\
STV    & 41.8958          & \textbf{0.9845}          & 31.7535          & \textbf{0.9285}          & 27.9748          & 0.8617          & 25.8979          & 0.8113          \\
WSTV   & \textbf{41.9860} & \textbf{0.9845} & \textbf{32.0720} & 0.9282 & \textbf{28.4314} & \textbf{0.8739} & \textbf{26.3218} & \textbf{0.8246} \\
Image  &  \multicolumn{8}{l}{Peppers}                                                                                                              \\
Models & PSNR             & SSIM            & PSNR             & SSIM            & PSNR             & SSIM            & PSNR             & SSIM            \\
TV     & 25.5367          & 0.9366          & 24.6559          & 0.8172          & 23.7426          & 0.7391          & 23.2832          & 0.7377          \\
ATV    & 40.6151          & 0.9706          & 31.8390          & 0.8838          & 28.2698          & 0.8142          & 26.1667          & 0.7562          \\
STV    & 41.7389          & 0.9782          & 32.5765          & \textbf{0.9067}          & 29.0054          & 0.8378          & 27.0236          & 0.7885          \\
WSTV   & \textbf{41.7927} & \textbf{0.9783} & \textbf{32.7879} & 0.9064 & \textbf{29.4400} & \textbf{0.8510} & \textbf{27.4317} & \textbf{0.8026} \\
Image  &  \multicolumn{8}{l}{Parrot}                                                                                                              \\
Models & PSNR             & SSIM            & PSNR             & SSIM            & PSNR             & SSIM            & PSNR             & SSIM            \\
TV     & 40.9081          & 0.9701          & 30.7923          & 0.8376          & 27.2613          & 0.7483          & 25.6409          & 0.7454          \\
ATV    & 40.4617          & 0.9695          & 30.8887          & 0.8769          & 27.4686          & 0.8040          & 25.6610          & 0.7551          \\
STV    & 41.7878          & 0.9811          & 31.3773          & \textbf{0.8951}          & 27.9851          & 0.8196          & 26.0868          & 0.7698          \\
WSTV   & \textbf{41.8670} & \textbf{0.9812} & \textbf{31.6405} & 0.8948 & \textbf{28.2974} & \textbf{0.8291} & \textbf{26.4374} & \textbf{0.7813} \\
Image  &  \multicolumn{8}{l}{House}                                                                                                              \\
Models & PSNR             & SSIM            & PSNR             & SSIM            & PSNR             & SSIM            & PSNR             & SSIM            \\
TV     & 40.3755          & 0.9612          & 31.8646          & 0.8035          & 28.8632          & 0.7185          & 27.9519          & 0.7420          \\
ATV    & 40.6282          & 0.9586          & 33.1480          & 0.8571          & 30.2274          & 0.8098          & 28.2457          & 0.7698          \\
STV    & 42.1200          & \textbf{0.9723}  & 33.7020         & \textbf{0.8726}  & 30.5140          & 0.8100          & 28.6823          & 0.7701          \\
WSTV   & \textbf{42.1399} & \textbf{0.9723} & \textbf{33.7861} & 0.8712 & \textbf{30.9519} & \textbf{0.8245} & \textbf{29.0902} & \textbf{0.7872} \\
Image  &  \multicolumn{8}{l}{Starfish}                                                                                                              \\
Models & PSNR             & SSIM            & PSNR             & SSIM            & PSNR             & SSIM            & PSNR             & SSIM            \\
TV     & 39.3999          & 0.9775          & 30.2435          & 0.8689          & 26.7045          & 0.7742          & 24.9694          & 0.7245          \\
ATV    & 40.2862          & 0.9792          & 30.3286          & 0.8820          & 26.5618          & 0.7833          & 24.4564          & 0.7024          \\
STV    & 41.4311          & \textbf{0.9854}          & 31.1023          & 0.9009          & 27.4800          & 0.8163          & 25.5208          & 0.7501          \\
WSTV   & \textbf{41.4639} & \textbf{0.9854} & \textbf{31.2817} & \textbf{0.9024} & \textbf{27.6146} & \textbf{0.8211} & \textbf{25.6366} & \textbf{0.7558} \\
\bottomrule\bottomrule
\end{tabular}}
\label{tab1}
\end{table}
			
\begin{table}[H]
\centering
\caption{Compare the PSNR and SSIM of different models in color image denoising.}
\resizebox{\textwidth}{!}{
\begin{tabular}{lllllllll}
\toprule\toprule
Noise  & 0.01    &       & 0.05    &         & 0.10    &         &0.15&                       \\ \midrule[1pt]
Image  &  \multicolumn{8}{l}{Airplane}                                                                                                              \\
Models & PSNR             & SSIM            & PSNR             & SSIM            & PSNR             & SSIM            & PSNR             & SSIM            \\
VTV    & 41.7862          & \textbf{0.9839}          & 31.5210          & 0.9185          & 27.8202          & 0.8622          & 25.8642          & 0.8006          \\
STV    & 41.8572          & 0.9828          & 31.7607          & 0.9159          & 27.9789          & 0.8406          & 26.0484          & 0.8053          \\
WSTV   & \textbf{42.0402} & 0.9838 & \textbf{32.1803} & \textbf{0.9274} & \textbf{28.4510} & \textbf{0.8733} & \textbf{26.4042} & \textbf{0.8237} \\
Image  &  \multicolumn{8}{l}{Baboon}                                                                                                              \\
Models & PSNR             & SSIM            & PSNR             & SSIM            & PSNR             & SSIM            & PSNR             & SSIM            \\
VTV    & 39.6139          & 0.9958          & 27.5152          & 0.9306          & 23.6604          & 0.8287          & 22.0993          & 0.7595          \\
STV    & 39.9461          & \textbf{0.9961}          & 27.8504          & 0.9358          & \textbf{24.1083}          & \textbf{0.8481} & 22.2328          & \textbf{0.7638} \\
WSTV   & \textbf{40.0153} & \textbf{0.9961} & \textbf{27.9325} & \textbf{0.9363} & 23.9745 & 0.8395          & \textbf{22.2351} & 0.7625          \\
Image  &  \multicolumn{8}{l}{Babara}                                                                                                              \\
Models & PSNR             & SSIM            & PSNR             & SSIM            & PSNR             & SSIM            & PSNR             & SSIM            \\
VTV    & 40.8893          & 0.9931          & 30.3346          & 0.9399          & 26.7817          & 0.8877          & 25.0124          & 0.8460          \\
STV    & 41.1515          & \textbf{0.9935} & 30.7391          & 0.9429          & 27.1864          & 0.8918          & 25.3847          & 0.8539          \\
WSTV   & \textbf{41.2126} & \textbf{0.9935} & \textbf{30.8368} & \textbf{0.9435} & \textbf{27.3048} & \textbf{0.8955} & \textbf{25.4745} & \textbf{0.8568} \\
Image  &  \multicolumn{8}{l}{Boats}                                                                                                              \\
Models & PSNR             & SSIM            & PSNR             & SSIM            & PSNR             & SSIM            & PSNR             & SSIM            \\
VTV    & 41.4372          & 0.9912          & 30.7406          & 0.9409          & 27.0070          & 0.8854          & 25.1392          & 0.8342          \\
STV    & 41.6649          & 0.9909          & 31.1500          & 0.9428          & 27.3744          & 0.8830          & 25.4324          & 0.8440          \\
WSTV   & \textbf{41.8123} & \textbf{0.9915} & \textbf{31.3970} & \textbf{0.9484} & \textbf{27.6057} & \textbf{0.8971} & \textbf{25.6314} & \textbf{0.8523} \\
Image  &  \multicolumn{8}{l}{Colored butterfly}                                                                                                              \\
Models & PSNR             & SSIM            & PSNR             & SSIM            & PSNR             & SSIM            & PSNR             & SSIM            \\
VTV    & 41.5569          & \textbf{0.9981}          & 31.0870          & 0.9831          & 27.2465          & 0.9628          & 25.0895          & 0.9413          \\
STV    & 41.6655          & 0.9980          & 31.2932          & 0.9835          & 27.3596          & 0.9622          & 25.2968          & 0.9436          \\
WSTV   & \textbf{41.8142} & \textbf{0.9981} & \textbf{31.6081} & \textbf{0.9848} & \textbf{27.7958} & \textbf{0.9663} & \textbf{25.6597} & \textbf{0.9474} \\
Image  &  \multicolumn{8}{l}{Colored peppers}                                                                                                              \\
Models & PSNR             & SSIM            & PSNR             & SSIM            & PSNR             & SSIM            & PSNR             & SSIM            \\
VTV    & 41.4440          & 0.9984          & 31.8107          & 0.9869          & 28.3800          & 0.9728          & 26.4494          & 0.9587          \\
STV    & 41.6158 & \textbf{0.9985} & 32.1747          & 0.9878          & 28.5737          & 0.9735          & 26.7581          & 0.9614          \\
WSTV   & \textbf{41.7445} & \textbf{0.9985} & \textbf{32.5572} & \textbf{0.9888} & \textbf{29.1003} & \textbf{0.9765} & \textbf{27.1139} & \textbf{0.9641} \\
\bottomrule\bottomrule
\end{tabular}}
\label{tab2}
\end{table}

\section{Conclusion}\label{Conclusion}
			
In this paper, we propose an image denoising model based on weighted structure tensor total variation. The core idea is that we use the anisotropic weighted matrix to the STV model to  characterize local features of the images,  which can effectively capture information about the restored images. In order to solve the corresponding WSTV model, we implement a fast first-order gradient projection algorithm for the dual optimization problem. The experimental results demonstrate that our method has  significant improvements compared to the TV-based and STV methods.
			
Although the WSTV model performs well, it takes a relatively long time when compared with other methods. Therefore, it is meaningful for us to study  the projection operators in the iteration process.  Customizing efficient algorithms according to a specific constraint set $\mathcal{C}$ and the parameter $q$ in the unit norm ball ${B}_{\mathcal{S}_{q}}$ will help WSTV reduce the running time.

\section* {Acknowledgments}
This research is supported by the National Natural Science Foundation of China (Grant No.  62073087).

\vspace{2ex}\noindent\textbf{Xiuhan Sheng}  received his BS degree in science from Henan Agricultural University in 2021. Currently, he is pursuing a master's degree at the School of Mathematics and Statistics, Guangdong University of Technology, Guangzhou, China. His research interest is mainly in image restoration.

\vspace{2ex}\noindent\textbf{Lijuan Yang}  received her BS degree in management from Hainan Medical University in 2021. Currently, she is pursuing a master's degree at the School of Mathematics and Statistics, Guangdong University of Technology, Guangzhou, China. Her research interest is mainly in machine learning.

\vspace{2ex}\noindent\textbf{Jingya Chang}  received her BS degree in science from Henan Normal University, her MS degree in science from Nanjing Normal University, and her PhD in applied mathematics from Hong Kong Polytechnic University in 2006, 2009, and 2018, respectively. Her research interests include nonlinear optimization, tensor optimization, and their applications in image processing and machine learning.

\listoffigures
\listoftables

\end{spacing}
\end{document}